\documentclass[onefignum,onetabnum]{siamart171218}

\usepackage{lipsum}
\usepackage{amsfonts}
\usepackage{graphicx}
\usepackage{epstopdf}
\usepackage{algorithmic}
\ifpdf
  \DeclareGraphicsExtensions{.eps,.pdf,.png,.jpg}
\else
  \DeclareGraphicsExtensions{.eps}
\fi

\newsiamremark{remark}{Remark}
\newsiamremark{hypothesis}{Hypothesis}
\crefname{hypothesis}{Hypothesis}{Hypotheses}
\newsiamthm{claim}{Claim}

\headers{Explicit Solutions and Geometric Algorithm}{H.-F Liu}

\title{Geometric Algorithm of Schr\"odinger Flow on a Sphere%\thanks{Submitted to the editors DATE.
}

\author{Hsiao-Fan Liu\thanks{Department of Mathematics, National Tsing Hua University, Taiwan 
  (\email{hfliu@math.nthu.edu.tw}).}}
\usepackage{amsopn}
\DeclareMathOperator{\diag}{diag}

\ifpdf
\hypersetup{
  pdftitle={Geometric Algorithm of Schr\"odinger Flow on a Sphere},
  pdfauthor={H.-F. Liu}
}
\fi

\begin{document}
\def\calL{\mathcal{L}}
\def\calG{\mathcal{G}}
\def\calD{\mathcal{D}}
\def\calJ{\mathcal{J}}
\def\calM{\mathcal{M}}
\def\calN{\mathcal{N}}
\def\calO{\mathcal{O}}
\def\calA{\mathcal{A}}
\def\calS{\mathcal{S}}
\def\calP{\mathcal{P}}
\def\calU{\mathcal{U}}
\def\calK{\mathcal{K}}
\def\frakgl{\mathfrak{gl}}
\def\frako{\mathfrak{o}}
\def\fraku{\mathfrak{u}}
\def\frakg{\mathfrak{g}}
\def\frakso{\mathfrak{so}}
\def\fraksl{\mathfrak{sl}}
\def\fraksp{\mathfrak{sp}}
\def\fraksu{\mathfrak{su}}
\def\F{\mathbb{F}}
\def\R{\mathbb{R}}
\def\N{\mathbb{N}}
\def\C{\mathbb{C}}
\def\M{\mathbb{M}}
\def\H{\mathbb{H}}
\def\P{\mathbb{P}}
\def\al{\alpha}
\def\be{\beta}
\def\p{\partial}
\def\n{\, | \, }
\def\ti{\tilde}
\def\a{\alpha}
\def\l{\lambda}
\def\hcalG {\hat{\mathcal{G}}}
\def\diag{{\rm diag \/ }}
\def\det{{\rm det \/ }}
\def\sp{{\rm span \/ }}
\def\rd{{\rm d\/}}
\def\K{\nabla}
\def\g{\gamma}
\def\Re{{\rm Re\/}}
\def\a{\alpha}
\def\b{\beta}
\def\d{\delta}
\def\D{\triangle}
\def\e{\epsilon}
\def\g{\gamma}
\def\G{\Gamma}
\def\K{\nabla}
\def\l{\lambda}
\def\L{\Lambda}
\def\n{\,\vert\,}
\def\o{\theta}
\def\w{\omega}
\def\W{\Omega}
\def\ca{{\mathcal{A}}}
\def\cd{{\mathcal{D}}}
\def\cf{{\mathcal{F}}}
\def\cg{{\mathcal{G}}}
\def\ch{{\mathcal{H}}}
\def\ck{{\mathcal{K}}}
\def\cl{{\mathcal{L}}}
\def\cL{{\mathcal{L}}}
\def\cm{{\mathcal{M}}}
\def\cn{{\mathcal{N}}}
\def\co{{\mathcal{O}}}
\def\cp{{\mathcal{P}}}
\def\cs{{\mathcal{S}}}
\def\ct{{\mathcal{T}}}
\def\cu{{\mathcal{U}}}
\def\cv{{\mathcal{V}}}
\def\cx{{\mathcal{X}}}
\def\li{\langle}
\def\ri{\rangle}
\def\n{\ \vert\ }
\def\tr{{\rm tr}}
\def\bs{\bigskip}
\def\ms{\medskip}
\def\ss{\smallskip}
\def\hb{\hfil\break\vskip -12pt}

\def\di{$\diamond$}
\def\ni{\noindent}
\def\ti{\tilde}
\def\p{\partial}
\def\Re{{\rm Re\/}}
\def\Im{{\rm Im\/}}
\def\I{{\rm I\/}}
\def\II{{\rm II\/}}
\def\diag{{\rm diag}}
\def\ad{{\rm ad}}
\def\Ad{{\rm Ad}}
\def\Iso{{\rm Iso}}
\def\Gr{{\rm Gr}}
\def\sgn{{\rm sgn}}

\def\rd{{\rm d\/}}

\def\R{\mathbb{R} }
\def\C{\mathbb{C}}
\def\H{\mathbb{H}}
\def\N{\mathbb{N}}
\def\Z{\mathbb{Z}}
\def\O{\mathbb{O}}
\def\F{\mathbb{F}}

\def\fg{\mathfrak{G}}

\newcommand{\beg}{\begin{example}}
\newcommand{\eeg}{\end{example}}
\newcommand{\bthm}{\begin{theorem}}
\newcommand{\ethm}{\end{theorem}}
\newcommand{\bprop}{\begin{proposition}}
\newcommand{\eprop}{\end{proposition}}
\newcommand{\bcor}{\begin{corollary}}
\newcommand{\ecor}{\end{corollary}}
\newcommand{\blem}{\begin{lemma}}
\newcommand{\elem}{\end{lemma}}
\newcommand{\bca}{\begin{cases}}
\newcommand{\eca}{\end{cases}}
\newcommand{\brem}{\begin{remark}}
\newcommand{\erem}{\end{remark}}
\newcommand{\bpm}{\begin{pmatrix}}
\newcommand{\epm}{\end{pmatrix}}
\newcommand{\bbm}{\begin{bmatrix}}
\newcommand{\ebm}{\end{bmatrix}}
\newcommand{\bvm}{\begin{vmatrix}}
\newcommand{\evm}{\end{vmatrix}}
\newcommand{\bdefn}{\begin{definition}}
\newcommand{\edefn}{\end{definition}}
\newcommand{\bsub}{\begin{subtitle}}
\newcommand{\esub}{\end{subtitle}}
\newcommand{\bex}{\begin{example}}
\newcommand{\eex}{\end{example}}
\newcommand{\ben}{\begin{enumerate}}
\newcommand{\een}{\end{enumerate}}

\newcommand{\balign}{\begin{align}}
\newcommand{\ealign}{\end{align}}
\newcommand{\baligns}{\begin{align*}}
\newcommand{\ealigns}{\end{align*}}
\newcommand{\beq}{\begin{equation}}
\newcommand{\eeq}{\end{equation}}
\newcommand{\beqa}{\begin{eqnarray}}
\newcommand{\eeqa}{\end{eqnarray}}

\def\pdo{$\psi$do}

\def\calA{{\mathcal A}}
\def\calB{{\mathcal B}}
\def\calF{{\mathcal F}}
\def\calG{{\mathcal G}}
\def\calJ{{\mathcal J}}
\def\calK{{\mathcal K}}
\def\calL{{\mathcal L}}
\def\calM{{\mathcal M}}
\def\calN{{\mathcal N}}
\def\calO{{\mathcal O}}
\def\calP{{\mathcal P}}
\def\calR{{\mathcal R}}
\def\calS{{\mathcal S}}
\def\calU{{\mathcal U}}
\def\calV{{\mathcal V}}

\def\li{\langle}
\def\ri{\rangle}

\def\frakP{{\mathfrak{P}}}

\def\half{\frac{1}{2}}
\def\Tr{{\rm Tr\/}}
\def\nkdv{$n\times n$ KdV}

\def \a {\alpha}
\def \b {\beta}
\def \d {\delta}
\def \D {\triangle}
\def \e {\epsilon}
\def \g {\gamma}
\def \G {\Gamma}
\def \K {\nabla}
\def \l {\lambda}
\def \L {\Lambda}
\def \n {\,\vert\,}
\def \N {\,\Vert\,}
\def \o {\theta}
\def\w{\omega}
\def\W{\Omega}
\def \s {\sigma}
\def \S {\Sigma}

\def\ca{{\mathcal {A}}}
\def\cC{{\mathcal {C}}}
\def\cg{{\mathcal {G}}}
\def\ci{{\mathcal {I}}}
\def\ck{{\mathcal {K}}}
\def\cl{{\mathcal {L}}}
\def\cm{{\mathcal {M}}}
\def\cn{{\mathcal {N}}}
\def\co{{\mathcal {O}}}
\def\cp{{\mathcal {P}}}
\def\cs{{\mathcal {S}}}
\def\ct{{\mathcal {T}}}
\def\cu{{\mathcal {U}}}

\def\R{{\mathbb{R}}}
\def\C{{\mathbb{C}}}
\def\H{{\mathbb{H}}}
\def\Z{{\mathbb{Z}}}

\def\Re{{\rm Re\/}}
\def\Im{{\rm Im\/}}
\def\tr{{\rm tr\/}}
\def\Id{{\rm Id\/}}
\def\I{{\rm I\/}}
\def\II{{\rm II\/}}
\def\li{\leftrangle}
\def\ri{rightrangle}
\def\id{{\rm Id}}
\def\gk{\frac{G}{K}}

\def\p{\partial}
\def\li{\langle}
\def\ri{\rangle}
\def\ti{\tilde}
\def\i{\/ \rm i }
\def\j{\/ \rm j }
\def\k{\/ \rm k}
\def\n {\ \vert\ }
\def\bu{$\bullet$}
\def\ni{\noindent}
\def\ii{{\rm i\,}}

\def\bs{\bigskip}
\def\ms{\medskip}
\def\ss{\smallskip}

\maketitle

% REQUIRED
\begin{abstract}
We construct the solution to the periodic Cauchy problem of the Schr\"odinger flow on the sphere. Such construction of solutions is formulated explicitly and therefore a geometric algorithm of solving this periodic Cauchy problem follows. Theoretical and experimental results will be discussed. 
\end{abstract}

% REQUIRED

% REQUIRED

\section{Introduction}

The Schr\"odinger flow on the unit sphere $\mathbb{S}^2$ is given by 
\beq\label{sf}
\gamma_t  =  \gamma \times \gamma_{xx},
\eeq
where $\g(x,t)$ is a real-valued vector function of $\mathbb{S}^2$ on $(x,t) \in \R \times [0,\infty)$, and $\times$ is the cross product in $\R^3$. We call $\g(x,t)$ satisfying \eqref{sf} a {\it Schr\"odinger curve}. This curve evolution \eqref{sf} corresponds to the non-linear Schr\"odinger equation (NLS)
\beq\label{nls} 
q_t= i (q_{xx} + 2 |q|^2 q).
\eeq
Takhtajan \cite{Ta77} first applied the inverse scattering method to describe its solution scheme. Tjon and Wright in \cite{TW77} studied solitons of \eqref{sf} in one dimension for the isotropic and anisotropic cases. The relation between \eqref{sf} and \eqref{nls} has been discussed by Terng and Uhlenbeck in \cite{TU06}, and  by Zakharov and Takhtadzhyan in \cite{ZT79}, where \eqref{sf} is known as the isotropic Heisenberg ferromagnet. 

It is well-known that \eqref{nls} is relevant to the vortex filament equation (VFE) in $\R^3$ defined by
\beq\label{vfe}
\a_t = \a_x \times \a_{xx}.
\eeq

This was first discovered by Da Rios in 1906 to model the movement of a thin vortex in a viscous fluid. In particular, \eqref{vfe} preserves the arc length parameter. So we may assume $||\a_x(x,t)||=1$ for all $t$, i.e., $\a(\cdot,t)$ is parametrized by its arc-length. Given a solution $\a(x,t)$ of VFE, Hasimoto (\cite{Has71} ) relates  \eqref{vfe} to \eqref{nls} via the following transform
\begin{equation}\label{za}
q(x,t)=k(x,t)e^{i(\theta(t)+\int_0^x \tau(s,t)ds)},
\end{equation}
where $k(\cdot, t),\tau(\cdot, t)$ are the curvature and torsion for $\a(\cdot, t)$, respectively, and $x$ is the arc-length parameter. Moreover, $q(x,t)$ is a solution of NLS. For example, a solution of the VFE $\a(x,t)= (\cos x, \sin x, t)$ with
$k(x,t)=1,$ $\tau(x,t)=0$ corresponds to a solution of NLS $q(x,t)= e^{2it}$.

Numerical solutions of VFE have been studied for years. In 1998, Hou, Klapper, and Si provided a formulation in \cite{HKS98} for calculating solutions of VFE numerically. This method is a generalization of their 2-D work \cite{HLS94}, in which they use the tangent angle ($\o$) and arc length ($L$) method for planar curves. This $\o-L$ formulation has none of high order time step stability constraints. Hou \emph{et al.} then proposed to use the normal principal curvatures $k_1, k_2$ as new variables to compute the motion of the curve in space. 

If one differentiates a solution $\a(x,t)$ of \eqref{vfe} with respect to $x$, a Schr\"odinger curve is given by $\a_x$. In general, numerical solution of \eqref{sf} can be obtained this way. Using geometry, Terng and Uhlenbeck in \cite{Ter15,TU06} formulate the correspondences between NLS, VFE and Schr\"odinger flow on $\mathbb{S}^2$. In fact, they systematically construct solutions of such curve flows by making use of Lax pairs of NLS and Lie theory.

This construction of solutions to a geometric curve flow, if found, leads to a geometric algorithm to obtain numerical solutions to a curve flow. Namely, this approach gives rise to both analytic and numerical solutions of Schr\"odinger flow directly. One advantage of this geometric method is that we reduce curve PDEs, which are highly nonlinear, to integrable partial differential equations and ODE systems. As one of the most famous integrable PDEs, the NLS can be solved and computed numerically using various methods (for example, the pseudo spectral method, which provides a good accuracy for periodic solutions). There are several schemes that solve ODE systems well, such as ode solvers in Matlab and Runge-Kutta method.

Another advantage is that there are  conservation properties since NLS is solitary. For instance, suppose $\g$ is a solution to \eqref{sf}.
We now consider the energy structure. Let $\mathcal{E}(\g(t))=|| \g_x(t)||^2_{L^2}$ be the energy function. Then the energy is conserved:
\beq\label{ener}
\mathcal{E}(\g(t))=\mathcal{E}(\g(0)),
\eeq
for any $t > 0$. Besides, there are infinitely many conservation laws for NLS. Hence, one can verify the accuracy of the algorithm by analyzing these conserved quantities. Finally, this geometric algorithm can be implemented easily, even for beginners in programming.

The paper is organized as follows.\

In \cref{sec:equi}, we recall the equivalence relations from results of Terng, Uhlenbeck and Thorbergsson. Proofs are given since they give rise to initial data in implementing. The essential idea works for both the Schr\"odinger flow and the VFE, although their constructions of formulations for curve solutions are different. We will focus on the Schr\"odinger flow on the unit sphere in the present article.
 In \cref{sec:main}, we formulate the solution of the periodic Cauchy problem of \eqref{sf} with an initial closed curve given. In \cref{sec:alg}, we give steps of our geometric algorithm, issues coming from implementing, error estimates for a stationary solution, and exhibit two explicit examples of periodic solutions with the viviani's curve and spherical sinusoid as initial curves. In the end of this section, we also give one example of VFE that features the motion of smoke rings. The discussions about error estimates are in \cref{sec:dis}, and the conclusions follow in
\cref{sec:conclusions}.

\section{Equivalence of the Schr\"odinger flow on Hermitian symmetric spaces and the NLS}
\label{sec:equi}

Suppose $(M, J, g, w)$ is a compact K\"ahler manifold with a complex structure $J$, the Riemannian metric $g$, and a symplectic form $w$ on $M$ satisfying $w(X,Y)= g(JX, Y)$. The {\it Schr\"{o}dinger flow} on $M$ (cf. \cite{TerTho01}) is the evolution equation
on $C^\infty(\R,M)$:
\beq\label{sch}
\gamma_t= J_{\gamma}(\nabla_{\gamma_x}\gamma_x),
\eeq
where $\nabla$ is the Levi-Civita connection of the metric $g$.
When $M=\C^n$, \eqref{sch} is the linear Schr\"odinger equation $\g_t= i\g_{xx}$. When $M=\mathbb S^2$, \eqref{sch} gives us the Heisenberg ferromagnetic model for
$\gamma: \R^2 \rightarrow \mathbb{S}^2$. Indeed, the complex structure of $\mathbb{S}^2$ at $\g$ sends $v$  to  $\gamma \times v$ on $ T\mathbb{S}^2_\gamma$, where $\times$ is the cross product in $\R^3$. Then
\beq
\gamma \times \nabla_{\gamma_{x}} \gamma_{x}=\gamma \times \gamma_{xx}^T=\gamma \times(\gamma_{xx}-(\gamma_{xx},\gamma)\gamma)=\gamma \times
\gamma_{xx},
\eeq 
which obviously is the evolution \eqref{sf} on $\mathbb{S}^2$.

The Schr\"{o}dinger flow \eqref{sch} is a Hamiltonian equation for the energy functional on $C^\infty(\mathbb{S}^1,M)$ with respect to an induced symplectic form by $\omega$ on $C^\infty(\mathbb{S}^1,M)$ (\cite{TerTho01}). Note that the critical points of the energy functional are geodesics of $(M,g)$, so the stationary solutions of the Schr\"odinger flow on $M$ are closed geodesics of $M$.
Furthermore, if $M$ is a Hermitian symmetric space, then \eqref{sch} can be written in terms of the Lie bracket. 

To be more precise, let $G$ be a simple complex Lie group, and
$\tau$ the involution that gives the maximal compact subgroup $U$. It is known that there
exists $a \in \calU$ such that $\ad(a)^2|_{\calP} = -\Id_{\calP}$ and $\calU = \calK \oplus \calP$, where $\calK$ is the centralizer of $a$ in $\calU$ and $\calP$ is the orthogonal complement of $\calK$. Then the Adjoint
$U$-orbit at $a$ in $\calU$ is diffeomorphic to $U$ and is a compact irreducible Hermitian $\frac{U}{K}$
symmetric space.

We briefly review some results on the Schr\"odinger flow proved by Terng and Uhlenbeck in \cite{TU06} for $\frac{U}{K}= \Gr(k, \C^n)$ and by Terng and Thorbergsson in \cite{TerTho01} for the other three classical Hermitian symmetric spaces. 
\begin{proposition} [\cite{TerTho01},\cite{TU06}] \label{la} Under the embedding of  the Hermitian symmetric space $\frac{U}{K}$ as the Adjoint orbit $U\cdot a$ in $\calU$, the Schr\"{o}dinger flow on $\frac{U}{K}$ is 
 \beq\label{ab}\g_t= [\g, \g_{xx}].
 \eeq 
 \end{proposition}
 
A Lax pair for \eqref{ab} is also derived and proved to be gauge equivalent to a Lax pair of NLS by Terng and Uhlenbeck \cite{TU06}. \cref{la} implies that if $M$ is a Hermitian symmetric space, then two evolutions, \eqref{sch} and \eqref{ab} are equivalent. In particular, we have $\g_t=[\g,\g_{xx}]=\g \times \g_{xx}$ for $M=\mathbb{S}^2$.

\cref{ac} and \cref{ac1} show the construction of NLS solutions from solutions of the Schr\"odinger flow and vice versa. We give a proof since the approach will be used later in coding.
\bthm[\cite{TerTho01},\cite{TerUhl00}]\label{ac} 
Let $\g:\R^2\to \frac{U}{K}$ be a solution of the Schr\"odinger flow on the Hermitian symmetric space $\frac{U}{K}=U\cdot a\subset \calU$. Then there exists $g:\R^2\to U$ satisfying
\ben
\item[(i)] $\g= gag^{-1}$,
\item[(ii)] $u= g^{-1}g_x:\R^2\to \calU_a^\bot$ satisfies the $\frac{U}{K}$-NLS equation:
\beq\label{gnls}
 u_t= [a, u_{xx}] -\frac{1}{2} [u, [u, [a, u]]],
 \eeq
 \item[(iii)]  $g^{-1}g_t= [a, u_x]-\frac{1}{2} [u, [a, u]]$. \een
 Moreover, $\ti g$ satisfies (i) and (ii) if and only if there is a constant $C \in U_a$ such that $\ti g=gC.$ 
\ethm
\begin{proof}
We recall that $K=U_a, P=U_a^\bot$ and $\calU=\calK \oplus \calP$.
Suppose $\g(x, t)$ is a solution of \eqref{ab}. Then
there exists $h: \R^2 \to U$ such that $\g(x, t) = h(x, t)ah(x, t)^{-1}$. Let $\pi_0, \pi_1$ be orthogonal projections of $\calU$ onto $\calK,\calP$, respectively. We choose $k: \R^2 \to K$ such that
$k_xk^{-1}=-\pi_0(h^{-1}h_x)$. Set $f(x, t)=h(x, t)k(x, t)$,
then $\gamma=faf^{-1}$. Moreover,
\begin{equation}
f^{-1}f_x=(hk)^{-1}(hk)_x=k^{-1}\pi_1(h^{-1}h_x)k \in
\calP.
\end{equation}
A direct computation shows that
$$
\gamma_x=f[f^{-1}f_x,a]f^{-1}=f[u,a]f^{-1} \mbox{ and } [\gamma
,\gamma_x]=f u f^{-1}.
$$
Since $\tau_\l=\gamma \l dx+(\gamma \l^2+[\gamma, \gamma_x]\l)dt$ is
flat for all $\lambda \in \C$, $f*\tau_\l$ is flat, i.e. the
following connection is flat for all $\lambda \in \C$:
\begin{equation}
f^{-1}\tau_\l f +f^{-1}df=(a\l+u)dx +(a\l^2+u\l+f^{-1}f_t)dt.
\end{equation}
Therefore, $(a\l^2 + u \l +f^{-1}f_t)_x-(a\l+u)_t+[a \l +u,
f^{-1}f_t]=0$.
\begin{align}\label{ld}
&u_x+[a,f^{-1}f_t]=0 \\
&(f^{-1}f_t)_x-u_t+[u,f^{-1}f_t]=0.
\end{align}

Write
\begin{equation*}
f^{-1}f_t=P+T,
\end{equation*}
where $P \in \calP$ and $T \in \calK$,
respectively. From \eqref{ld}, we have
\begin{equation}
P=[a,u_x], ~T_x= -\frac{1}{2}[u,[a,u]]_x.
\end{equation}
So, $T= -\frac{1}{2}[u,[a,u]]+c(t)$ for some function $c(t)$.

Define $g=fy(t)$, where $y(t) \in \calK$ such that
$y_ty^{-1}=-c(t).$

Next, we will show that $g$ defined above satisfies the conditions $(i)-(iii)$. Since $y(t)$ and $a$ commute, it is easy to see that
$gag^{-1}=\gamma$.  In particular,
\begin{align}
&g^{-1}g_x=y^{-1}f^{-1}f_xy=y^{-1}uy \in \calP, \\
&g^{-1}g_t=y^{-1}(f^{-1}f_t+y_ty^{-1})y=-\frac{1}{2}[y^{-1}uy,[a,y^{-1}uy]],
\end{align}
which means $y^{-1}uy$ is a solution of \eqref{gnls}.

For the uniqueness, suppose $\tilde{g}$ satisfies $(1)-(2)$, and set
$C=g^{-1}\tilde{g}$. Then
\beq
\tilde{g}^{-1}\tilde{g}_x=C^{-1}g^{-1}g_xC+C^{-1}C_x.
\eeq
Since $\tilde{g}^{-1}\tilde{g}_x$ and $C^{-1}g^{-1}g_xC$ are in
$\calP$ while $C^{-1}C_x \in \calK$,
\beq
C^{-1}C_x=0.
\eeq
Similarly, $C^{-1}C_t=0$. So $C$ is constant.

\end{proof}

The converse is also true.
\bthm[\cite{TerTho01}, \cite{TerUhl00}]\label{ac1}  
Let $u:\R^2\to \calU_a^\bot$ be a smooth solution of \eqref{gnls}. Then given any $c_0\in U$,
the following linear system for $g:\R^2\to U$,
\beq
\bca
 g^{-1}g_x= u,\\
 g^{-1}g_t= [a, u_x] - \frac{1}{2}[u,[a,u]], \\
 g(0,0)= c_0
\eca
\eeq
has a unique smooth solution $g:\R^2 \to U$. Moreover,  $\g(x,t)= g(x,t) a g(x,t)^{-1}$ is a solution of the Schr\"odinger flow \eqref{ab} on $\frac{U}{K}$. 
\ethm
In fact, when $\l=\l_0$ is any arbitrary real number, a shift of $\g=gag^{-1}$ by $2\l_0$ is also a solution of \eqref{ab}.
\bprop\label{ac2}
Let $u: \R^2 \to \calU_a^\bot$ be a solution of \eqref{gnls} and $E$ satisfy
\beq\label{saa}
\bca
 E^{-1}E_x=a\l+ u,\\
 E^{-1}E_t=a\l^2+u\l+ [a, u_x] - \frac{1}{2}[u,[a,u]].
\eca
\eeq
If $\l_0 \in \R$ and $g(x,t)=E(x,t,\l_0)$, then $\g=gag^{-1}(x-2\l_0t,t)$ is a solution of \eqref{ab}.
\eprop
\begin{proof}
Let $\eta(x,t)=gag^{-1}(x,t)$ and denote $Q_{-1}=[a, u_x] - \frac{1}{2}[u,[a,u]]$.
It can be checked that
\beq
\begin{array}{l}
\eta_x=g[u,a]g^{-1},\\
 \eta_t=g[u\l_0+Q_{-1},a]g^{-1}.
\end{array}
\eeq
Direct computations show that
\beq
\gamma_{xx}=g[a\l_0+u,[u,a]]g^{-1}+g[u_x,a]g^{-1},
\eeq
and therefore we obtain
\beq
\gamma \times \gamma_{xx} = g[a,u\l_0]g^{-1}+gu_xg^{-1}.
\eeq
We see that $\gamma_t=-2\l_0\eta_x+\eta_t$, which gives
\beq
g[-u\l_0,a]g^{-1}+g[Q_{-1},a]g^{-1}.
\eeq
Here, since $[Q_{-1},a]=u_x$, $\gamma_t=[\gamma, \gamma_{xx}].$
\end{proof}
We say $E$ is a frame of solution $u$ if $E$ is a solution of the ODE system \eqref{saa}.  
\brem
Consider $\g(x,t)$ to be a great circle on the unit sphere, i.e., $\g(x,t)=(\cos x,\sin x,0)$. It is obviously a solution to \eqref{ab}. Following from the proof of Theorem \ref{ac}, one local invariant $q(x,t)$ is given by
\beq
q(x,t)=\frac{i}{2}e^{\frac{i}{2}t},
\eeq
that of course can be checked to solve the NLS.
\erem

Soliton solutions of the Schr\"odinger equation have been widely studied in history. One can compute soliton solutions using explicit formulations that are known. In fact, B\"acklund transformation gives $N$-soliton solutions for integrable systems, especially for the NLS. For example, \cite{LTW} says the following:
\bthm\label{ac22}(B\"acklund transformation)
Let $E(x,t,\l)$ be a frame of a solution $u$ of the NLS \eqref{nls}, $\a \in \C\setminus \R$, and $\pi$ a Hermitian projection of $\C^2$ onto $V$.
 Define 
 \beq
 k_{\a,\pi}(\l)=\I_n+\frac{\a-\bar \a}{\l-\a}\pi^\perp.
 \eeq
 Let $\ti \pi(x,t)$ be the Hermitian projection of $\C^2$ onto $\ti V(x,t)= E(x,t,\a)^{-1}V$.   Then we have 
  \beq
  \ti u= u+ (\a-\bar \a)[a, \ti\pi]
  \eeq
  is a solution of the NLS and 
  \beq\ti E(x,t,\l)= k_{\a, \pi}(\l) E(x,t,\l) k^{-1}_{\a, \ti\pi(x,t)}\eeq is a frame of $\ti u$. 
\ethm

\section{Main results}
\label{sec:main}
The correspondence discussed in \cref{ac} and \cref{ac1} gives a systematic way to construct solution of the Schr\"odinger flow, and hence the initial value problem of the Schr\"odinger flow can be solved. In this section, we focus on the sphere case, in particular, we consider the periodic Cauchy problem. Namely, given a closed curve $\g_0(x)$ on a unit sphere, we show the existence of a $L$-periodic curve that evolves according to the Schr\"odinger flow \cref{sf} and write down the algebraic solution formula explicitly. Without loss of generality, we assume the period $L$ is $2\pi$ for the rest of article. The main theorem is the following.

\bthm\label{zu}
Given a smooth closed curve $\gamma_0(x) : [0, 2 \pi] \rightarrow \mathbb{S}^2$ with period $2\pi$. Then there exists a unique $x$-periodic $\gamma(x,t)$ with period $2 \pi$ satisfying

\begin{equation}\label{ivp1}
\left\{
\begin{array}{rcl}
\gamma_t & = & \gamma \times \gamma_{xx}\\
\gamma(x,0)& =&  \gamma_0(x)
\end{array}
\right..
\end{equation}
\ethm

We first note that Theorem \ref{zu} is a special case of the following Theorems \ref{ac} and \ref{ac1}. 
Theroem \ref{ac} shows that, for any arbitrary $\g_0(x): [0, 2 \pi] \to \mathbb{S}^2$ given, there is a {\it frame} $f:\R \to SU(2)$ such that $\gamma_0 = f a f^{-1}$, where $a=\diag(\frac{i}{2},-\frac{i}{2})$, satisfying $f(0)=I_2$ and $f^{-1}f_x = u_0$, where $u_0$ is of the form
\beq
u_0 = 
\left(
\begin{array}{cc}
0&q_0\\
-\bar{q_0}&0
\end{array}
\right).
\eeq
We notice that $f(x)$ found may not be periodic, so the essential idea is to find a periodic one. Below we show how to construct a periodic frame of $\g_0$. 

Since $\gamma_0$ is periodic, $\gamma_0(2\pi)=\gamma_0(0)$. It yields that $f(2\pi)a=af(2\pi).$
That is, $f(2\pi)$ lies in the centralizer $SU(2)_a=\{\diag(e^{\ii\o},e^{-\ii \o}) \n \o \in [0,2\pi)\}$ and hence we may write 
\beq
f(2\pi) = e^{2\pi c_0a},
\eeq
for some constant $c_0$. A direct computation gives the following proposition.
\bprop\label{ea}
Define 
\beq
\tilde{f}(x) = f(x)e^{-c_0ax}.
\eeq
Then $\tilde{f}(x)$ has the following properties:
\begin{enumerate}
\item $\gamma_0=\tilde{f}a \tilde{f}^{-1}$
\item $\tilde{f}(x)$ is periodic in $x$
\item $\tilde{f}^{-1}\tilde{f}_x =
         \left(
         \begin{array}{cc}
         -\frac{i}{2}c_0 & \tilde{q_0}\\
         -\bar{\tilde{q_0}} & \frac{i}{2}c_0
         \end{array}
         \right)
         $, where $\tilde{q_0} (x)= q_0(x)e^{ic_0x}.$
         \item $\tilde{q_0}$ is periodic.
\end{enumerate}
\end{proposition}
\cref{ea} gives us a way to decompose the periodic curve $\g_0$ in terms of a periodic {\it frame} $\ti f(x)$ and a local invariant $\ti q_0(x)$. Here, we call $c_0$ the {\it normal holonormy}. Next, we evolve the curve according to the partial differential equation 
\beq\label{gpde}
\g_t=\g \times \g_{xx}.
\eeq
\bprop
Suppose $\g(x,t):[0,2\pi] \to \mathbb S^2$ solves \eqref{gpde} and is periodic in $x$ with periodic $2\pi$. By Theorem \ref{ac}, there exists $f:\R^2 \to SU(2)$ such that $\g=faf^{-1}$,  $f^{-1}f_x=u$, and $f^{-1}f_t=Q_{-1}$, where 
\beq
a=\diag(\frac{i}{2},-\frac{i}{2}), \quad u=\bpm0&q\\-\bar q&0\epm, \quad Q_{-1}=\frac{i}{2}\bpm-|q|^2&q_x\\\bar q_x&|q|^2\epm.
\eeq
Define $c_0(t)$ to be a function of $t$ satisfying 
\beq\label{ea1}f^{-1}(0,t)f(2\pi,t)=e^{2\pi c_0(t) a}.\eeq
Then $c_0(t)$ is independent of $t$.
\eprop
\begin{proof}
Taking $t$-derivative of \eqref{ea1} gives 
\begin{align*}
e^{2\pi c_0(t)a}2\pi c_0'(t)a&=-f^{-1}(0,t)f_t(0,t)f(2\pi,t)+f^{-1}(0,t)f_t(2\pi,t)\\
&=e^{2\pi c_0(t)a}Q_{-1}(2\pi,t)-Q_{-1}(0,t)e^{2\pi c_0(t)a}.
\end{align*}
It is easy to see that 
\beq
2\pi c_0'(t)a=Q_{-1}(2\pi,t)-e^{-2\pi c_0(t)a}Q_{-1}(0,t)e^{2\pi c_0(t)a}.
\eeq
A direct computation shows 
\beq
e^{-2\pi c_0(t)a}Q_{-1}(0,t)e^{2\pi c_0(t)a}=\frac{i}{2}\bpm-|q|^2&q_xe^{-4\pi i c_0(t)}\\ \bar q_x e^{4\pi i c_0(t)}& |q|^2\epm.
\eeq
Note that $Q_{-1}(0,t)=Q_{-1}(2\pi,t)$. So, $c_0'(t)=0,$ as desired.
\end{proof}

Next, we consider the periodic Cauchy problem for NLS. Suppose that $q:\R^2 \to \C$ is a solution \footnote{This periodic Cauchy problem has a global smooth solution from the results in \cite{B93,Bo94,Its76}).} of 
 \begin{equation}\label{eb}
\left\{
\begin{array}{ccl}
q_t & = & i(q_{xx}+2|q|^2q)\\
q(x,0)& = & \tilde{q_0}(x).
\end{array}
\right.,
\end{equation}

Let $E$ satisfy
\beq\label{ey}
\bca
E^{-1}E_x=\bpm\frac{i}{2}\l & q\\-\bar q& -\frac{i}{2} \l \epm,\\
E^{-1}E_t = \bpm\frac{i}{2}\l^2-i |q|^2 & q\l+iq_x \\-\bar q\l+i\bar{q}_x & -\frac{i}{2} \l^2+i|q|^2 \epm,\\
E(0,0,\bar \l)^*=E(0,0,\l)^{-1}.
\eca
\eeq
Then it turns out that there is a periodic frame for a solution $q$ of NLS periodic in $x$ in the following proof. 

%\bthm\label{pc}
%Given a smooth and periodic curve $\gamma_0 : [0, 2 \pi] \rightarrow \mathbb{S}^2$ with $\g_0(0)=a$. Then there exists a unique $\gamma(x,t)$ periodic in $x$ with period $2 \pi$ satisfying \eqref{ivp1}.
%\ethm
\begin{proof}[Proof of \cref{zu}]
Without loss of generality, we assume $\g_0(0)=a$. We know that there is $f \in SU(2)$ such that $\g_0 = f a f^{-1}$ and 
\beq
f^{-1}f_x=\bpm0& q_0\\-\bar q_0 &0 \epm.
\eeq 
Since $\g_0$ is periodic, $f(2\pi)$ commutes with $a$. So 
$f(2\pi)=e^{2\pi c_0 a}$ for some $c_0\in \R.$ Define 
\beq
\ti f(x)=f(x)e^{-c_0 ax}.
\eeq
By Proposition \ref{ea}, $\ti f$ is periodic and $\g_0=\ti f a \ti f^{-1}$. In particular, 
\beq
\ti f^{-1}\ti f_x = \bpm -\frac{i}{2}c_0& q_0(x)e^{i c_0x}\\ -\bar q_0(x)e^{-i c_0x}& \frac{i}{2}c_0  \epm.
\eeq

Let $q(x,t)$ be the solution of 
\beq
\left\{
\begin{array}{ccl}
q_t & = & i(q_{xx}+2|q|^2q)\\
q(x,0)& = & q_0(x)e^{ic_0x}
\end{array}
\right.,
\eeq
periodic in $x$, and $E(x,t,\l)$ the extended frame for $q$ satisfying
\beq\label{ivpE}
\bca
E^{-1}E_x=a (-c_0)  + u,\\
E^{-1}E_t = a c_0^2 +u (-c_0) + Q_{-1}(u),\\
E(0,0,-c_0) =  \tilde{f}(0).
\eca
\eeq
We claim that $g(x,t) = E(x,t,-c_0)$ is periodic in $x$ with period $2\pi$.
Let $y(t)= g(2\pi,t)-g(0,t)$. We know $g^{-1}g_t = c_0^2a-c_0u+Q_{-1}(u)$ and $u=\bpm 0& q\\-\bar q &0\epm$ is periodic. Then
\begin{align*}
y'(t)&=g(2\pi,t)(c_0^2a-c_0u+Q_{-1}(u))|_{x=2\pi}-g(0,t)(c_0^2a-c_0u+Q_{-1}(u))|_{x=0}\\
&=(g(2\pi,t)-g(0,t))(c_0^2a-c_0u+Q_{-1}(u))|_{x=0}\\
&=y(t)A(t),
\end{align*}
where $A(t)=(c_0^2a-c_0u+Q_{-1}(u))|_{x=0}.$

Since $y(0)=0$ solves the ODE $y'(t)=y(t)A(t)$, the uniqueness theorem of ODE shows that $y(t) \equiv 0$. The claim follows.
 Let $\eta=gag^{-1}$. Then $\g(x,t)=\eta(x+2c_0t,t)$ is a solution of $\g_t=\g \times \g_{xx}$ by \cref{ac2}. 

It remains to verify the initial condition. Note that Proposition \ref{ea} implies
\beq
\gamma(x,0) = \eta(x,0)=\tilde{f}(x)a\tilde{f}^{-1}(x)=\gamma_0(x).
\eeq
In particular, that $\gamma$ is periodic in $x$ follows from the periodicity of $E(x,t,-c_0)$. Finally, the uniqueness of $\gamma$ follows from the uniqueness of $E(x,t,-c_0)$.
\end{proof}

\section{Algorithm and Experimental results}
\label{sec:alg}

Our analysis leads to the following algorithm to solve numerically the periodic Cauchy problem \eqref{zu} with initial data $\g_0(x):[0,2\pi] \to \mathbb S^2$ being a closed initial curve on $\mathbb S^2$. In summary, the programming steps are as follows.
\ben
\item[ $1$.] We first write $\g_0$ as an element in $\fraksu(2)$ and diagonalize $\g_0$ to find $f \in SU(2)$ such that $\g_0=faf^{-1}$ and $$f^{-1}f_x=\bpm0&q_0\\-\bar q_0&0\epm.$$
\item[ $2$.] Then compute $c_0$ by solving $f^{-1}(0)f(2\pi)=e^{2\pi c_0a}$. 
\item[ $3$.] We use the WGMS method in \cite{RSP} (implicit spectral method) to solve the periodic Cauchy problem of NLS \eqref{eb} with the initial data $q_0(x)e^{ic_0x}.$
\item[ $4$.] compute $E$ by solving the ODE system \eqref{ivpE} with the right hand side given by solutions $q$ of \eqref{eb} and the initial data $fe^{-c_0ax}$.
\item[ $5$.] calculate $\g=EaE^{-1}$ in terms of elements in $\fraksu(2)$ and then we map them back to $\R^3$, which is the numerical solution to \eqref{zu}.
\een

Next, we compute the normal holonormy $c_0$, which gives the initial periodic data $q_0(x)e^{ic_0 x}$ for the periodic Cauchy problem of NLS. The pseudo-spectral method \cite{RSP} provides the numerical periodic solution to NLS with great accuracy. 

With such numerical local invariant $q$ fitting in the right hand side of \eqref{ey}, we solve \eqref{ey} using the Runge-Kutta method. Once $E$ is obtained, we calculate $\g=EaE^{-1}$ in terms of elements in $\fraksu(2)$ and then we map them back to $\R^3$. By Proposition \ref{ac2} and the interpolation, the numerical solution to \eqref{zu} is derived.
\subsection{Experimental issues}
The first step can be done simply by diagonalizing $f$. One can easily check that the diagonal entries are $0.5i$ and $-0.5i$, however, the function {\bf eig}$(\cdot)$ in Matlab does not work well here. On one hand, the analytic $f$ in Step $1$ is a function that satisfies $\g_0(x)=f(x)af(x)^{-1}$ at each point $x$. On the other hand, we plug grid points $x_i$'s into $\g_0(x)$ before applying the {\bf eig} in Matlab. This makes the {\bf eig} function treat $\g_0(x_i)$'s as individual scalar matrices, and then returns eigenvectors of $\g_0(x_i)$ that don't obey the same function since matrices formed by eigenvectors are not unique. For this reason, we need another way to figure out the $f$.

Identifying $\R^3$ with the skew-Hermitian matrices $\fraksu(2)$, we see that $\fraksu(2)$ has a standard basis consisting of the four elements 
\beq
a=\bpm \frac{i}{2}&0\\0&-\frac{i}{2} \epm, b=\bpm 0&\frac{1}{2}\\-\frac{1}{2}&0\epm, \mbox{ and }c=\bpm 0&\frac{i}{2}\\ \frac{i}{2}&0\epm,
\eeq
and the map between elements in $\fraksu(2)$ and vectors $\g_0=(r_1,r_2,r_3)$ on the sphere is
\beq
\g_0=(r_1,r_2,r_3)\mapsto r_1a+r_2b+r_3c= \bpm \frac{i}{2}r_1&\frac{1}{2}r_2+\frac{i}{2}r_3\\ -\frac{1}{2}r_2+\frac{i}{2}r_3&-\frac{i}{2}r_1\epm,
\eeq
denoted by $\Gamma_0$.
A standard calculation of eigenvectors for $\Gamma_0$ shows
\beq F=
\bpm
\sqrt{\frac{1+r_1}{2}} & \frac{i}{\sqrt{2}}\frac{r_2+i r_3}{\sqrt{1+r_1}}\\
 \frac{i}{\sqrt{2}}\frac{r_2-i r_3}{\sqrt{1+r_1}} &\sqrt{\frac{1+r_1}{2}}
\epm
\eeq
It can be checked that $F \in SU(2)$ and such $F$ is not unique. We make such choice of $F$ for the following reasons. When $r_1=1$, $r_2$ and $r_3$ are obviously zero, i.e., $\Gamma_0(1)=a$. This immediately implies the matrix of eigenvectors is the $2 \times 2$ identity matrix, which agrees with our formulation of $F$. However, when $r_1=-1$, $\Gamma_0(-1)=-a$. The eigenvectors are 
\beq
\bpm 0& i\\ i& 0\epm.
\eeq
Besides, let $r_1\rightarrow -1$, we have
\beq
\frac{r_2+i r_3}{\sqrt{1+r_1}} \rightarrow \sqrt{2} e^{i\o}, \quad \o=\tan \frac{r_3}{r_2}.
\eeq
The limit of $F$ does not exist as $r_1\rightarrow -1$, and hence, this formulation is not continuous at $r_1=-1$. Without loss of generality, let $\g_0$  not pass through the point $(-1,0,0)$.
 
We also notice that $F^{-1}F_x$ might not be off-diagonal matrix. In order to make this happen, we follow Theorem \cref{ac} to rotate $F$ by a matrix $K$ where $-K_xK^{-1}$ is equal to the diagonal terms of $F^{-1}F_x$. For instance, consider $\g_0(x)=(0,\cos x, \sin x)$, then 
\beq
F(x) = \bpm \frac{1}{\sqrt{2}} e^{\frac{i}{2}x}&  \frac{i}{\sqrt{2}} e^{\frac{i}{2}x} \\
 \frac{i}{\sqrt{2}} e^{-\frac{i}{2}x}&  \frac{1}{\sqrt{2}} e^{-\frac{i}{2}x} \epm, \quad 
 F^{-1}F_x=\bpm 0&-\frac{1}{2}\\ \frac{1}{2} & 0 \epm.  
\eeq

\subsection{Experimental errors}

While actually implementing, the errors are provided to verify the accuracy of this geometric scheme. We would also like to see if numerical solutions obtained from our geometric algorithm preserve these relevant quantities. A fixed point $\g=a$ is the trivial solution of \eqref{sf} with the local invariant $q=0$. Our implement immediately shows a fixed point on the sphere if we start with $\g_0(x)=(1,0,0)$. Given the initial curve to be a circle, then the exact solution is stationary, i.e., \beq\g(x,t)=(0,\cos x, \sin x).\eeq The corresponding invariant is $q(x,t)=-\frac{1}{2}e^{\frac{i}{2}t}$, a solution of NLS. The following errors are estimated between the numerical solution $\g_{num}(x,t)$ and $\g(x,t)$ by the $L^2$-norm at each time $t$:\\
\beq
E_N(t) = \left(\int_{0}^{2\pi} |\g_{num}(x,t)-\g(x,t)|^2 ~dx\right)^{\frac{1}{2}},
\eeq
where $N$ indicates the number of $x$ partitions. \\

\begin{table}[htbp]
\caption{$\triangle t=0.001$}\label{table.1}
\begin{center}
\begin{tabular}{ccccc}
\hline\hline
$t$\\\mbox{ steps}        & $E_{64}(t)$        & $E_{128}(t)$   & $E_{256}(t)$ & $E_{512}(t)$    \\
\hline
10 & 1.0851962E-02 &	2.722708E-03 & 6.84783037E-04 &	1.72447709E-04\\
20 & 1.0257888E-02  &	2.582122E-03 &	6.51822917E-04 &	1.58102531E-04\\
30 & 9.711601E-03 &	2.454641E-03 &	6.14796552E-04 &	1.34178170E-04 \\
40 & 9.220668E-03  &	2.338847E-03 &	5.70308441E-04 &	1.19274466E-04 \\
50 & 8.79326E-03 &	2.234083E-03 &	5.19933864E-04 &	1.16633012E-04 \\
60 & 8.43756E-03  &	2.140637E-03 &	4.72659359E-04 &	1.19427212E-04 \\
70 & 8.160922E-03 &	2.059924E-03 &	4.47994913E-04 &	1.44297698E-04 \\
80 & 7.968884E-03  &	1.994621E-03 &	4.59273153E-04 &	1.93683390E-04  \\
90 & 7.864274E-03 &	1.94872E-03 &	4.94939577E-04 &	2.49089493E-04  \\
100 & 7.846684E-03  & 1.9274E-03 &	5.43976410E-04 &	3.11981115E-04  \\
\hline
\end{tabular}
\end{center}
\end{table}
We also consider the maximal difference between $\g_{num}(x,t)$ and $\g(x,t)$ at each time $t$ given by
\beq
E_N^{sup}(t) = max_{x \in I}\left(\|\g_{num}(x,t)-\g(x,t)\|\right),
\eeq
where $I=\{j\frac{2\pi}{N} | j=0,1,2, \cdots, N-1\}.$ 

\begin{table}[htbp]
\caption{$\triangle t=0.001$}\label{table.2}
\begin{center}
\begin{tabular}{ccccc}
\hline\hline
$t$\\\mbox{ steps}      & $E^{sup}_{64}(t)$        & $E^{sup}_{128}(t)$   & $E^{sup}_{256}(t)$ & $E^{sup}_{512}(t)$    \\
\hline
10 & 2.50304072E-04	&	8.83938616E-05	&	3.87637202E-05	&	1.59494744E-05	 \\
20 & 4.10069784E-04	&	1.60699569E-04	&	6.49229140E-05	&	1.97234863E-05	 \\
30 & 5.69783116E-04	&	2.22613774E-04	&	7.95435964E-05	&	1.68690680E-05	 \\
40 & 7.26480096E-04	&	2.74553194E-04	&	8.42922932E-05	&	2.76624306E-05	 \\
50 & 8.79882976E-04	&	3.17302096E-04	&	8.29884807E-05	&	4.21233398E-05	 \\
60 & 1.03000937E-03	&	3.51880564E-04	&	8.41160665E-05	&	5.42705369E-05	 \\
70 & 1.17668715E-03	&	3.79578845E-04	&	1.00586167E-04	&	7.04118124E-05	 \\
80 & 1.31944921E-03	&	4.02029934E-04	&	1.32082869E-04	&	9.19189874E-05	 \\
90 & 1.45756158E-03	&	4.21296895E-04	&	1.64970296E-04	&	1.13954497E-04	 \\
100& 1.59011707E-03	&	4.39944903E-04	&	1.95175345E-04	&	1.37739549E-04	 \\
\hline
\end{tabular}
\end{center}
\end{table}%
It is obvious that at a fixed time, the errors $E_N(t)$ and $E_N^{sup}(t)$ are decreasing when $N$ becomes bigger. Similarly,  it is natural for us to consider how $E_N(t)$ behaves in time when $N$ is fixed. In \cref{table.3}, we demonstrate the global errors for different sizes of time steps. Namely, 
\beq
G_N^{sup} = max_{0\leq t_n \leq T}\left(E_N^{sup}(t_n) \right),
\eeq
where $T$ is the total time period.
\begin{table}[htbp]
\caption{The error $G_{N}^{sup}$ computed for different $\triangle t$ with $T=1$}\label{table.3}
\begin{center}
\begin{tabular}{ccccc}
\hline\hline
$N$      & $\triangle t=\frac{1}{200}$        & $\triangle t=\frac{1}{400}$    & $\triangle t=\frac{1}{800}$     & $\triangle t=\frac{1}{1600}$\\
\hline
64	&	2.43337736E-03	&	1.94066115E-03	&	1.76206517E-03	&	1.63673643E-03\\
128	&	1.45758341E-03	&	8.77318416E-04	&	5.73463467E-04	&	4.49949019E-04\\
256	&	1.27041144E-03	&	6.73414782E-04	&	3.71644911E-04	&	2.23587772E-04\\
512	&	1.22892162E-03	&	6.29339213E-04	&	3.25308022E-04	&	1.69441660E-04\\
1024	&	1.21935787E-03	&	6.19703034E-04	&	3.15140360E-04	&	1.59348657E-04\\
\hline
\end{tabular}
\end{center}
\end{table}%

\begin{figure}[htbp]
  \centering
  \label{fig:a}\includegraphics[scale=0.5]{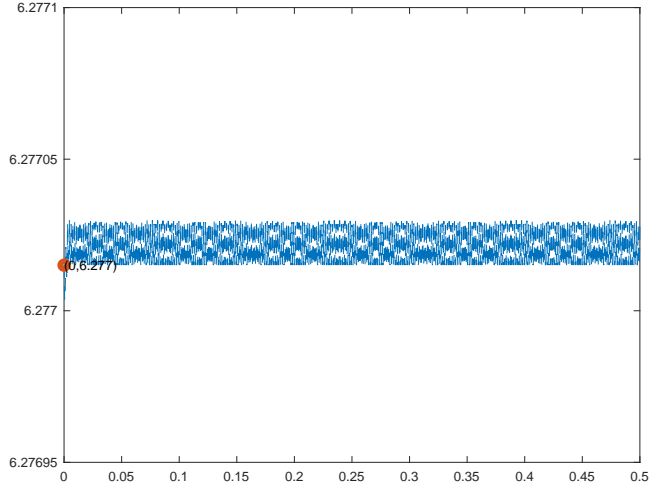}
  \caption{Numerical energy $\mathcal{E}(\g_{num}(t))$ at each time $t$ versus the initial energy $6.277$ at $t=0$ with $\triangle t =10^{-4}$ and total time $0.5$.}
\end{figure}

As we can see in \cref{fig:a}, the numerical energy at each time is approximately equal to the initial energy. The difference of energy at each time is less than $0.0005$ according to \cref{fig:a}. However, the real energy in this case is $2\pi$, which is differed from the numerical energy within $0.01$. The error occurs due to accumulation error of trapezoidal integration of numerical solutions $\g_{num}(x,t)$ and machine error. The experimental results do give us a numerical energy closer to $2\pi$ if $\triangle t, \triangle x$ become smaller.

It is well-known that the NLS has infinitely many conserved quantities. The first four  conserved quantities for NLS are:

\beq\displaystyle H_1=\oint |q|^2 ~dx, \quad H_2=\oint \bar{q}q_x ~dx,\eeq
\beq\displaystyle H_3=\oint |q_x|^2-|q|^4~dx, \quad H_4=\oint q\bar q_x-\bar q q_x~dx.\eeq
Although one can compute the errors for the conserved quantities, it is expected that the inaccuracy will raise while calculating the integral using the trapezoidal method and the derivatives of $q$.

\subsection*{Viviani's curve} 
\begin{figure}[htbp]
  \centering
  \label{fig:d}\includegraphics[scale=0.7]{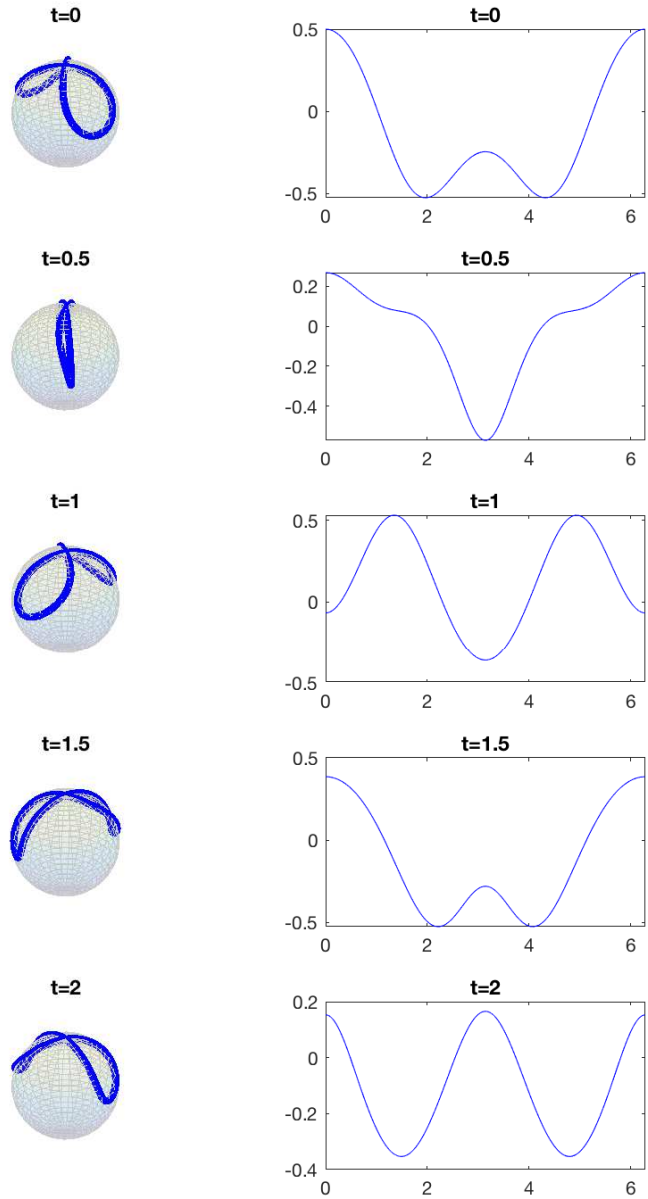}
  \caption{Numerical solution $\g_{num}$ versus the real part of the corresponding local invariant $q$ at $t=0, 0.5, 1, 1.5, 2$, respectively with $\triangle t=0.001$.}
\end{figure}

If we start with $$\g_0(x)=(\sin x \cos x,\sin x,\cos^2 x),$$ which has the figure eight shape. It is also considered to be the intersection of a sphere centered at the origin with a cylinder tangent to the sphere and passes through the origin. If one projects such curve stereographically from the point diametrically opposite the double point, then the lemniscate of Bernoulli is obtained. \cref{fig:d} gives the motion of Schr\"odinger curve from $t=0$ to $t=2$ with time step $0.001$ and $N=2^{10}$. The right column in \cref{fig:d} consists of the behavior of corresponding local invariant $q$.

\subsection*{Spherical Sinusoid} 
\begin{figure}[htbp]
  \centering
  \label{fig:i}\includegraphics[scale=0.8]{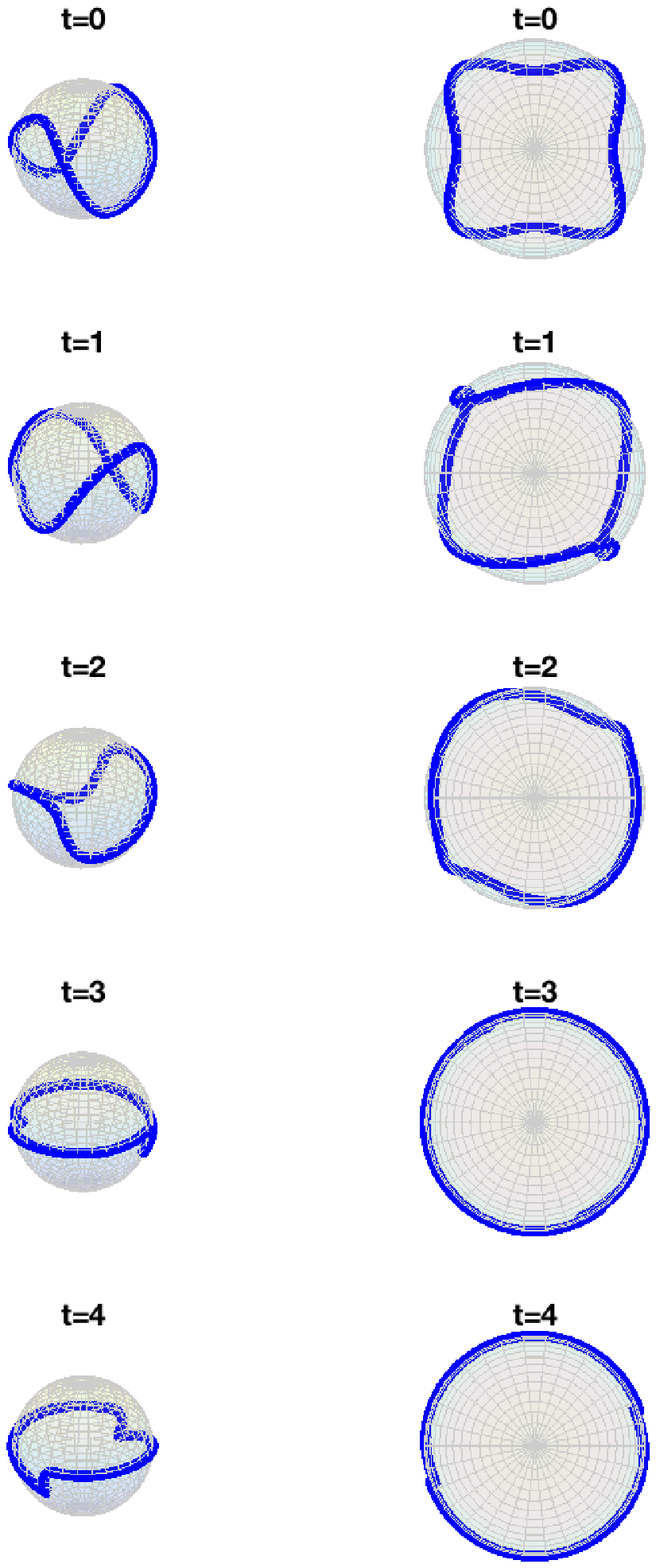}
  \caption{Numerical solution $\g_{num}$ (left) versus the projection of the $\g_{num}$ onto $xy$ plane (right) at $t=0,1,2,3,4$, respectively with $\triangle t=0.01$.}
\end{figure}

\begin{figure}[htbp]
  \centering
  \label{fig:j}\includegraphics[scale=0.65]{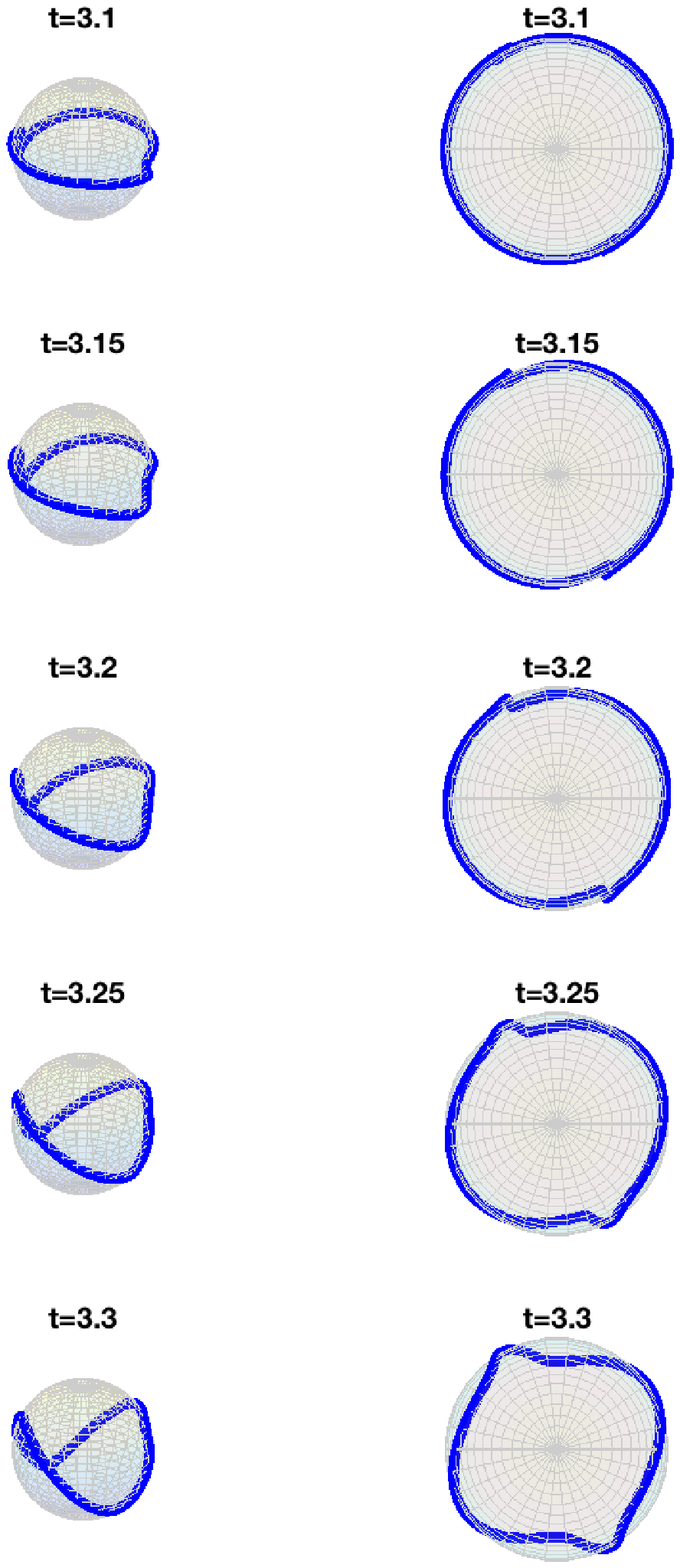}
  \caption{Numerical solution $\g_{num}$ (left) versus the projection of the $\g_{num}$ onto $xy$ plane (right) at $t=3.1,3.15,3.2,3,25,3.3$, respectively with $\triangle t=0.01$.}
\end{figure}

Next example is to begin with 
$$\g_0(x)=(\frac{\cos x}{\sqrt{1+\cos^2 2x}},\frac{\sin x}{\sqrt{1+\cos^2 2x}},\frac{\cos 2x}{\sqrt{1+\cos^2 2x}}).$$ 
\cref{fig:i} shows the numerical results. The left column consists of Schr\"odinger curves obtained from the initial curve $\g_0(x)$ with $N=2^{10}$ and the time step $0.01$ at different time, paired with those curves from overhead viewpoint. At $t=3$ and $t=4$, the curves seem to have cusps only because of different perspectives. They are actually smooth. In \cref{fig:j}, the drop-down parts of curve stretches. So we see that the drop-down parts slightly move to other places at each time in the right column of \cref{fig:j}.

\subsection*{Schr\"odinger curve with $1$-soliton $q$} 
\begin{figure}[htbp]
  \centering
  \label{fig:g}\includegraphics[scale=0.7]{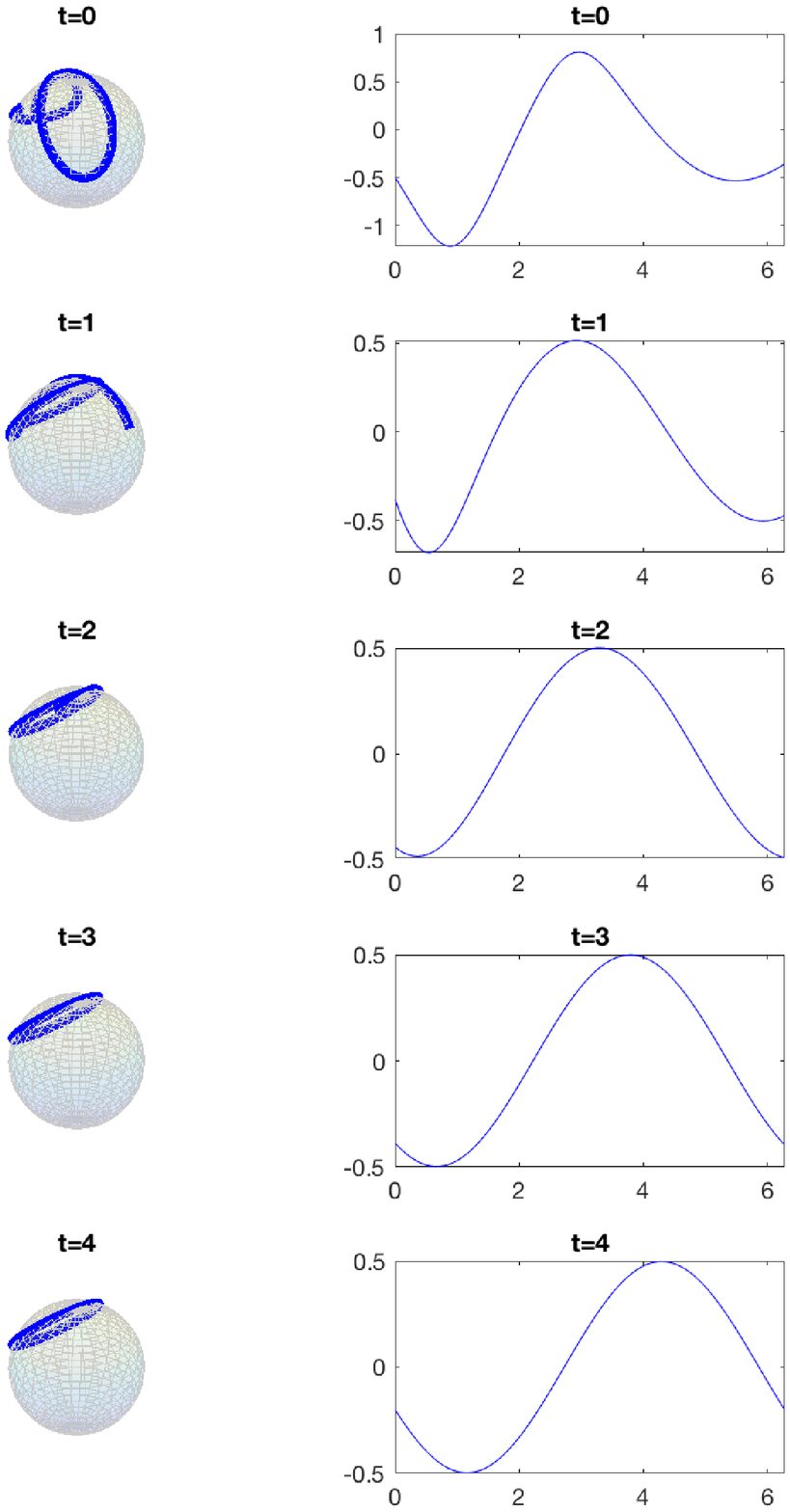}
  \caption{$\ti \g_{num}$ constructed from applying BT on the circle versus the real part of the corresponding local invariant $q$ at $t=0, 1, 2, 3, 4$, respectively with $\triangle t=0.01$.}
\end{figure}

Based on our geometric scheme, we are able to obtain the first periodic solution $q$ of the NLS. Applying B\"acklund transformation to $q$ will give us one-soliton $\ti q$.
\cref{ac2} and \cref{ac22} imply that the new $\ti E(x,t,\l)$ will give rise to a new solution of \eqref{sf}. As an example, we apply BT on the stationary solution, i.e., a great circle. \cref{fig:g} shows a numerical result with $\a=1-i, V=\C\left(\begin{array}{c}1 \\i\end{array}\right)$.    

\begin{figure}[htbp]
  \centering
  \label{fig:h}\includegraphics[scale=0.5]{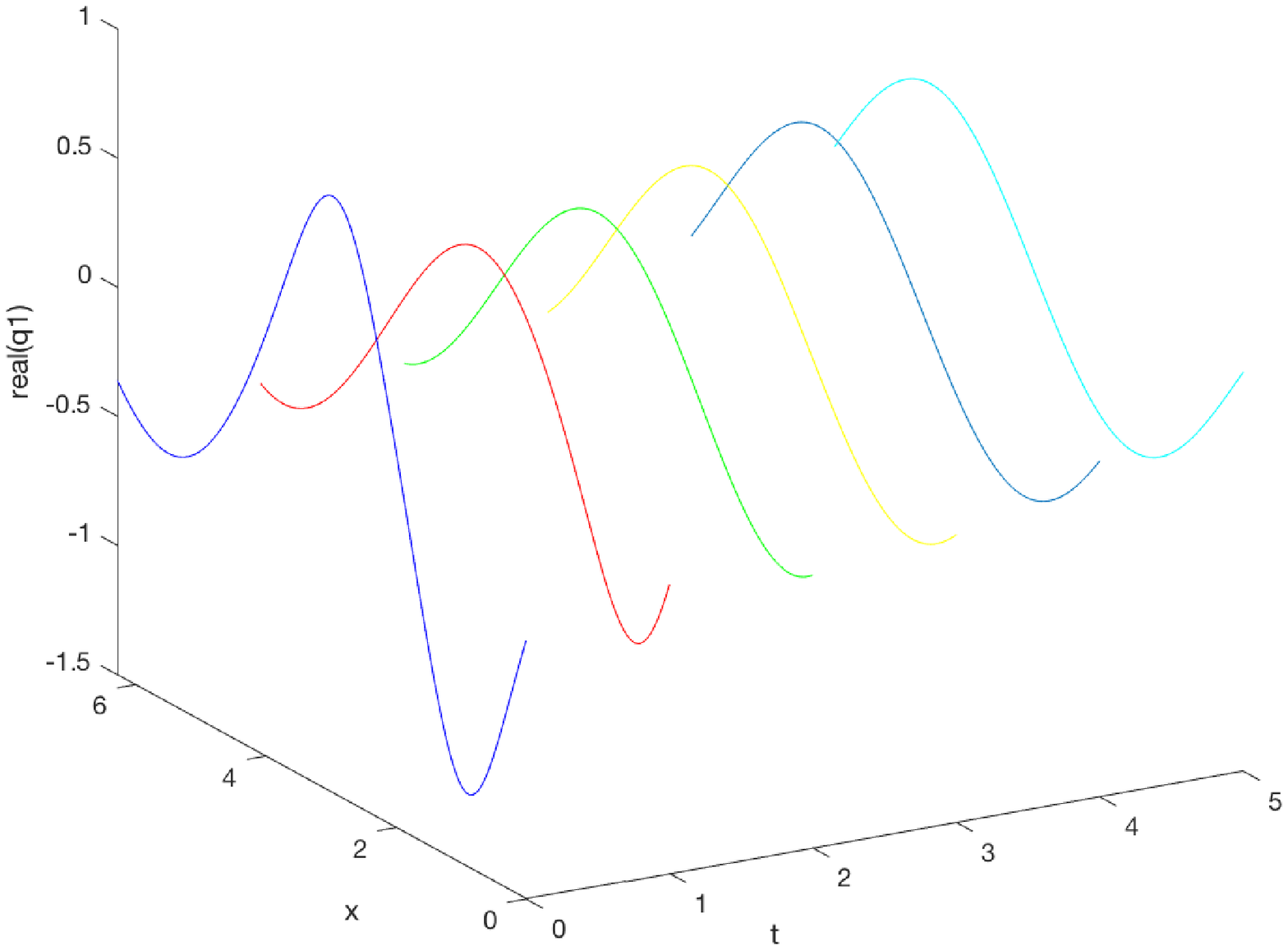}
  \caption{The real part of $1$-soliton solution $q_1$ corresponding $\ti \g_{num}$  constructed from the circle with $q$ at $t=0, 1, 2, 3, 4, 5$, respectively with $\triangle t=0.01$.}
\end{figure} 
\brem

The first picture in \cref{fig:g} indicates that the solution we get by applying the B\"acklund transformation might not be periodic. At $t=1$, two endpoints are approaching to each other. From the experimental outcomes, such curve remains the same circle as that on the last picture in \cref{fig:g} until $t=10$.
\erem

\subsection{Vortex filament equation} 
Assume $\g(x,t)$ is a solution of the Schr\"odinger flow \eqref{sf}. Let 
\beq\label{za2}
\a(x,t) =\int_0^x \g(s,t)~ds+c(t),
\eeq
where $c(t)$ is independent of $x$. It implies
\begin{align*}
\a_t&=\int_0^x \g_t(s,t)~ds +c'(t)\\
&=\int_0^x \g \times \g_{xx}(s,t)~ds+c'(t).
\end{align*}
Note that the integrand is a total derivative, i.e., 
\beq\g \times \g_{xx} = (\g \times \g_x)_x,\eeq
so the above equation turns out to be
\beq\label{za1}
\a_t=\g \times \g_x(x,t)-\g \times \g_x(0,t)+c'(t).
\eeq
It says from \eqref{za1} that $\a(x,t)$ is a solution of \eqref{vfe} if and only if 
\beq
c'(t) =\g \times \g_x(0,t).
\eeq
Moreover, what we can say for periodicity is the following:
\bprop
$\a(x,t)$ defined as \eqref{za2} is periodic in $x$ with period $2\pi$ if a solution of the Schr\"odinger flow $\g(x,t)$ has the same $x$-period.
\eprop
\begin{proof}
Let \beq y(t) = \a(2\pi,t)-\a(0,t).\eeq Direct computation shows that 
\beq y'(t) = \g \times \g_x(2\pi,t)-\g \times \g_x(0,t).\eeq
Since $\g(x,t)$ is periodic in $x$ with period $2\pi$, so is $\g_x(x,t)$. It is clearly that $y'(t)=0$, namely, $y(t)$ is a constant. As $y(0)=0$, we obtain $y(t)=0$ for all $t$. The assertion is proved. 
\end{proof}

An immediate example is when a great circle $\g(x,t)=(0,\cos x, \sin x)$ solves the Schr\"odinger flow, one can verify that $\a(x,t)=(t,\sin x, -\cos x)$ gives a solution of the vortex filament equation. Together with our geometric scheme mentioned in previous sections, a numerical solution of the vortex filament equation \eqref{vfe} has been provided.

On the other hand, an algebraic construction of periodic solution for the vortex filament equation
\beq\label{ee}
\bca
\a_t=\a_x \times \a_{xx}\\
\a(x,0)=\a_0
\eca,
\eeq
has been given by Terng where  $\a_0 : [0,2\pi] \to \R^3$ is a smooth arc-length parametrized curve periodic in $x$ with period $2\pi$. Equivalently, a geometric scheme follows from such construction in  \cite{Ter15}. We summarize her ideas without proof.
\bthm\label{es}
Let $\a(x,t)$ be a solution of the VFE \eqref{vfe} that is periodic in $x$ with period $2\pi$ and $\N \a_x \N =1$. Suppose $(e_0,\vec{n}_1,\vec{n}_2)$ is orthonormal along $\a$ such that $e_0=\a_x.$ Let $\w=(\vec{n}_1)_x\cdot \vec{n}_2$. Then 

\beq\label{ee3}
c_0=\frac{1}{2\pi} \int_0^{2\pi} \w(x,t)~dx
\eeq 
is constant for all $t$, and there is $g=(u_0,u_1,u_2)(x,t)$ such that 
\ben
\item $g(\cdot,t)$ is a periodic h-frame along $\a(\cdot,t)$,
\item $g^{-1}g_x=\left(\begin{array}{ccc}0 & - \zeta_1 & -\zeta_2 \\ \zeta_1 & 0 & -2c_0 \\ \zeta_2 & 2c_0 & 0\end{array}\right)$,
\item $q=\frac{1}{2}(\zeta_1 + i \zeta_2)$ solves the nonlinear Schr\"odinger equation 
\beq\label{ee2}
q_t=\frac{i}{2}(q_{xx}+2|q|^2q).
\eeq
\een
\ethm
\bprop\label{fa}
Let $q$ be a solution of \eqref{ee2} periodic in $x$ with period $2\pi$, $\l_0\in \R$, and $E(x,t,\l)$ the extended frame of $q$. If $E(x,0,\l_0)$ is periodic in $x$ with period $2\pi$, then so is $E(x,t,\l_0)$.
\eprop
\bprop\label{eu1}
Suppose $\a_0(x):\R^2 \to \R^3$ is a periodic curve parametrized by arc-length with period $2\pi$. Let $(u_0^0,u_1^0,u_2^0)$ be a $h$-frame $2q_1^0,2q_2^0$ and $\phi \in SU(2)$ such that \beq(u^0_0,u^0_1,u^0_2)=(\phi a \phi^{-1}, \phi b \phi^{-1}, \phi c  \phi^{-1}).\eeq Suppose $q:\R^2 \to \C$ is a periodic solution of 
\beq\label{}
\bca
q_t= \frac{i}{2}(q_{xx}+2 | q|^2q),\\
q(x,0)=q_1^0+i q_2^0.
\eca
\eeq
Let $E$ be the extended frame with initial data $E(0,0,c_0)=\phi, \eta =E_\l E^{-1}\n_{\l=c_0},$ and $
\a(x,t)=\eta(x-2c_0 t,t)$. Then $\tilde \a(x,t)= \a(x,t) -\eta(0,0) + \a_0(0)$ solves \eqref{ee} and is periodic in $x$ with period $2\pi$.
\eprop

The geometric algorithm to solve \eqref{ee} follows from the above discussion, we demonstrate experimental results here.
\begin{figure}[htbp]
  \centering
  \label{fig:f}\includegraphics[scale=0.35]{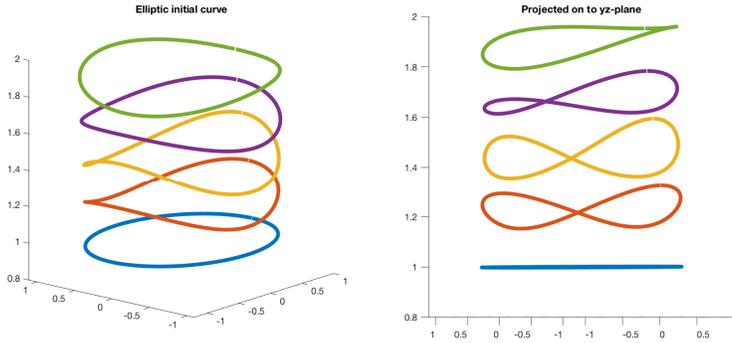}
  \caption{Curve motion starting with the bottom curve $\a_0(x)=(\cos x,\sin x,\cos x)$.}
\end{figure}
Starting with the bottom curve $\a_0(x)=(\cos x,\sin x,\cos x)$ (left), this curve slowly floats up as time goes by. Such simulation shows that it captures the {\it smoke ring} feature for the vortex filament equation.

\section{Discussion of error estimates}\label{sec:dis}
At a fixed time, \cref{table.1} shows that the $L^2$-error $E_N(t)$ decreases when the number of grid points $N$ is increase. Indeed, from the experimental results, we see that $\frac{E_{2N}(t)}{E_N(t)}$ is approximately located in between $\frac{1}{4}$ and $\frac{1}{2}$.  

We also notice that if $N$ is fixed, the error $E_N(t)$ accumulates at each time step in \cref{table.1}. For the total time $T=1$, we see the errors $G_N^{sup}$ in \cref{table.3}. As for the data demonstrated, the error becomes half of itself when we reduce time step by a half, provided $N$ is fixed. The error estimates shown are not so impressed in a sense of numerics. One reason why this error is not "so small" is that we use the WGMS method to approximate the solution $q$ to the NLS. Several types of numerical errors come from the WGMS algorithm, including the obvious error in approximating the integral and a truncation error when the fixed point iteration was stopped after a finite number of steps. 

Another reason is that we simply use finite difference method to obtain numerical results for the frame $E$ in \eqref{ivpE} with the right hand side filled out with the estimated $q$. This of course increases errors. It also indicates the numerical scheme can be improved by choosing other more accurate ODE solvers.

\section{Conclusions}
\label{sec:conclusions}

Although the accuracy provided is not relatively impressed, this geometric algorithm shows that using simple solvers for each piece in the implement can help to get numerical solutions "good enough" to the nonlinear curve PDE \eqref{sf}. The advantage of this method is transform the nonlinearity of curve motion to solving the ODE system \eqref{ivpE}. Each step can be solved numerically by built-in functions in MatLab, therefore our scheme is easier for beginners to do programming. The price to pay from the experimental results is obviously some accuracy. However, it is expected that experts in coding can have better approximations if they work more on each step of the implement.

%\appendix
%\section{An example appendix} 
%\lipsum[71]
%
%\begin{lemma}
%Test Lemma.
%\end{lemma}

\section*{Acknowledgments}
The author would like to thank Prof. Chuu-Lian Terng at University of California, Irvine, for her guidance on theoretical results, and Yu-Yu Liu at the Department of Mathematics, National Cheng Kung University for his useful discussion on numerical techniques.

\bibliographystyle{siamplain}
\bibliography{references}
\end{document}